\documentclass[10pt]{amsart}
\usepackage{amssymb,amsmath,amsfonts,amsthm,graphics,mathrsfs,amscd,hyperref}
\usepackage[hmargin=1in,vmargin=1in]{geometry}
\usepackage[all,cmtip]{xy}

\theoremstyle{plain}

\newtheorem{lemma}[equation]{Lemma}

\numberwithin{equation}{subsection}

\title{Poisson deformations of ruled surfaces over an elliptic curve}
\author{Chunghoon Kim}

\email{ckim042@gmail.com}            

\begin{document}

\maketitle
\begin{abstract}
We determine obstructedness or unobstructedness of (holomorphic) Poisson deformations of ruled surfaces over an elliptic curve.
\end{abstract}

\tableofcontents

We study (holomorphic) Poisson deformations of ruled surfaces over an elliptic curve. A holomorphic Poisson manifold $X$ is a complex manifold such that its structure sheaf is a sheaf of Poisson algebras.\footnote{For general information on Poisson geometry, we refer to \cite{Lau13} .}  A holomorphic Poisson structure is encoded in a holomorphic section (a holomorphic bivector field) $\Lambda_0\in H^0(M,\wedge^2 \Theta_M)$ with $[\Lambda_0,\Lambda_0]=0$, where $\Theta_M$ is the sheaf of germs of holomorphic vector fields on $M$, and the bracket $[-,-]$ is the Schouten bracket on $M$. In the sequel a holomorphic Poisson manifold will be denoted by $(M,\Lambda_0)$. In \cite{Kim15}, we studied deformations of holomorphic Poisson structures on the basis of Kodaira-Spencer's deformation theory of complex structures. We defined a concept of a family of compact holomorphic Poisson manifolds, called a Poisson analytic family, which is based on a complex analytic family in the sense of Kodaira-Spencer's deformations theory (\cite{Kod05}). Throughout this paper, we will call deformations of complex structures `complex deformations', and deformations of holomorphic Poisson structures `Poisson deformations' for short. Given a compact holomorphic Poisson manifold $(M,\Lambda_0)$, infinitesimal (Poisson) deformations of $(M,\Lambda_0)$ are encoded in the first cohomology group $\mathbb{H}^1(M, \Theta_M^\bullet)$ of the complex of sheaves $\Theta_M^\bullet: \Theta_M \to \wedge^2 \Theta_M\to\wedge^3 \Theta_M \to \cdots $ induced by $[\Lambda_0,-]$. We say that a compact holomorphic Poisson manifold $(M,\Lambda_0)$ is unobstructed in Poisson deformations if there is a Poisson analytic family $(\mathcal{M},\Lambda, B, \pi)$ of deformations of $\pi^{-1}(0)=(M,\Lambda_0), 0\in B$ such that the associated Poisson Kodaira-Spencer map $\varphi_0: T_0 B \to \mathbb{H}^1(M,\Theta_M^\bullet)$ is an isomorphism at $0\in B$. Otherwise, we say that $(M,\Lambda_0)$ is obstructed in Poisson deformations. In this paper we determine obstructedness or unobstructedness of Poisson deformations for ruled surfaces $S$ over an elliptic curve $X$ for any holomorphic Poisson structure on $S$. Our method is based on \cite{Suw69} in the context of Poisson deformations. 

In section $\ref{section1}$, we review a biholomorphic classification of ruled surfaces over an elliptic curve and explicit constructions of biholomorphic equivalent classes of ruled surfaces $S$ over an elliptic curve $X$ presented in \cite{Suw69}. We explicitly describe $H^0(S,\Theta_S)$ and $H^1(S,\Theta_S)$ in terms of a covering of $S$. In section $\ref{section2}$, based on the constructions in section $\ref{section1}$, we explicitly describe holomorphic Poisson structures $H^0(S,\wedge^2 \Theta_S)$ on $S$, and $H^1(S,\wedge^2 \Theta_S)$. Then by using them, in section $\ref{section3}$, we compute $\mathbb{H}^0(S,\Theta_S^\bullet), \mathbb{H}^1(S,\Theta_S^\bullet)$ and $\mathbb{H}^2(S,\Theta_S^\bullet)$ for any Poisson ruled surface $(S,\Lambda_0)$ over an elliptic curve $X$ (see Table $\ref{ruled}$). Finally in section $\ref{section4}$, we determine obstructedness or unobstructedness of Poisson deformations for any Poisson ruled surface $(S,\Lambda_0)$ over an elliptic curve $X$ (see Table $\ref{ruled}$). For unobstructed Poisson ruled surfaces $(S,\Lambda_0)$ over an elliptic curve, we explicitly construct Poisson analytic families of deformations of $(S,\Lambda_0)$ such that the associated Poisson Kodaira-Spencer map is an isomorphism at the distinguished point.

\section{Ruled surfaces over an elliptic curve}\label{section1}

 We review a biholomorphic classification of ruled surfaces over an elliptic curve presented in \cite{Suw69} to which we refer for the details. Every ruled surfaces over an elliptic curve $X$ can be expressed uniquely as one of the following: (i) a line bundle of non-negative degree, (ii) $A_0$, and (iii) $A_{-1}$ where $A_0$ and $A_{-1}$ are affine bundles (\cite{Ati55}, \cite{Ati57}). For any divisor $D$ on $X$, we denote by $[D]$ the line bundle over $X$ which is determined by $D$. We write the operation of the group of divisors on $X$ multiplicatively. Then by Abel-Jocobi theorem, given a point $p_0$ on $X$, the mapping $p\to [p_0 p^{-1}]$ gives an isomorphism between complex torus $X$ and the Picard variety $Pic^0(X)=ker\, c\cong H^1(X,\mathcal{O}_X)/hH^1(X,\mathbb{Z})$, where
\begin{align}\label{a1}
\cdots \to H^1(X,\mathbb{Z})\cong \mathbb{Z}\oplus \mathbb{Z}\xrightarrow{h} H^1(X, \mathcal{O}_X)\cong \mathbb{C}\xrightarrow{e} H^1(X,\mathcal{O}_X^*)\xrightarrow{c} H^2(X,\mathbb{Z})\cong \mathbb{Z}\to 0,
\end{align} 
which is the long exact sequence induced from the exponential sequence $0\to \mathbb{Z}\to \mathcal{O}_X\to \mathcal{O}_X^*\to 0$,  and for any line bundle $\xi$ of degree $n\geq 1$, there exists a point $p$ on $X$ such that $\xi=[p^n]$. It follows that all the ruled surfaces associated with line bundles of degree $n\geq 1$ are biholomorphically equivalent to one and the same surface, which will be denoted by $S_n$. We denote by $\mathcal{S}_0$ the biholomorphic equivalent classes of ruled surfaces associated with line bundles of degree $0$. Then every ruled surfaces over an elliptic curve $X$ can be classified biholomorphically as follows:
\begin{align}\label{a60}
\mathcal{S}_0,\,\,\,\,\,S_n(n\geq 1), \,\,\,\,\,A_0, \,\,\,\,\,A_{-1}.
\end{align}

Let $S$ belong to one of  $(\ref{a60})$. We will explicitly construct $S$ and describe $H^0(S,\Theta_S)$, and $H^1(S,\Theta_S)$ (for the details, see \cite{Suw69} Theorem $3$ and Theorem $4$). We represent an elliptic curve $X$ as a quotient group : $X=\mathbb{C}/G_\omega$, where $G_\omega$ is a discontinuous group of the additive group $\mathbb{C}$ generated by $\omega$ and $1$ with $\text{Im}\,\omega>0$, and for any $u\in \mathbb{C}$, we denote by $[u]$ the corresponding element of $X=\mathbb{C}/G_\omega$. Take a point $p\in X$ and let $u_1$ be a local coordinate at $p$. Let $U_1=\{u_1 | |u_1|<\epsilon\}$ for a sufficiently small number $\epsilon>0$ and $U=X-p$. Then $\mathscr{U}=\{U,U_1\}$ is the Stein covering of $X$. We construct $S$ by gluing $U\times \mathbb{P}_\mathbb{C}^1$ and $U_1\times \mathbb{P}_\mathbb{C}^1$ in the following way. Here we denote by $\xi$ and $\xi_1$ the inhomogeneous coordinates of $\mathbb{P}_\mathbb{C}^1$ and $\mathbb{P}_\mathbb{C}^1$ respectively in the covering $\,\,\mathcal{U}:=\{U\times \mathbb{P}_\mathbb{C}^1,U_1\times \mathbb{P}_\mathbb{C}^1\}$ of $S$. 

\subsection{Construction of $S=S_0=X\times \mathbb{P}_\mathbb{C}^1$ and descriptions of $H^0(S,\Theta_S)$ and $H^1(S,\Theta_S)$} \label{subsection1}\

We set $S=(U\times \mathbb{P}_\mathbb{C}^1)\cup (U_1\times \mathbb{P}_\mathbb{C}^1)$, where $(u,\xi)\in U\times \mathbb{P}_\mathbb{C}^1$ and $(u_1,\xi_1)\in  U_1\times \mathbb{P}_\mathbb{C}^1$ are identified if and only if $\xi=\xi_1$, and $[u]=p+u_1$. Then $\dim_\mathbb{C} H^0(S,\Theta_S)=4$ and
\begin{align}\label{a7}
\frac{\partial}{\partial u},\,\,\,\,\,\frac{\partial}{\partial \xi},\,\,\,\,\,\xi\frac{\partial}{\partial \xi},\,\,\,\,\,\xi^2\frac{\partial}{\partial \xi}\in H^0(S, \Theta_S)
\end{align}
forms a basis of $H^0(S, \Theta_S)$. On the other hand, $\dim_\mathbb{C}H^1(S,\Theta_S)=4$ and
\begin{align}\label{a8}
\frac{1}{u_1}\frac{\partial}{\partial u_1},\,\,\,\,\,\frac{1}{u_1}\frac{\partial}{\partial \xi_1},\,\,\,\,\,\frac{\xi_1}{u_1}\frac{\partial}{\partial \xi_1},\,\,\,\,\,\frac{\xi_1^2}{u_1}\frac{\partial}{\partial \xi_1}\in C^1(\mathcal{U}, \Theta_S)
\end{align}
forms a basis of $H^1(S,\Theta_S)$.
\subsection{Construction of $S\in \mathcal{S}_0, S\ne S_0$, and descriptions of $H^0(S,\Theta_S)$ and $H^1(S,\Theta_S)$}\label{subsection2}\ 

We note that any line bundle of degree zero can be represented by a $1$-cocycle $\eta(t)=\{\eta_{ij}(t)\}_{i,j=0,1}\in C^1(\mathscr{U},\mathcal{O}_X^*)$, $\eta_{01}(t)=e^{\frac{t}{u_1}}$ for some $t\in \mathbb{C}$, where $\mathscr{U}=\{U,U_1\}$. $\eta(t)$ represents the trivial bundle if and only if $\eta'(t)=\{\eta'_{ij}(t)\}_{i,j=0,1}\in C^1(\mathscr{U},\mathcal{O}_X)$, where $\eta'_{01}(t)=\frac{t}{2\pi i u_1}$ is in the image of $h$ in $(\ref{a1})$. In this case, we say that $t$ belongs to the lattice.

We set $S=(U\times \mathbb{P}_\mathbb{C}^1)\cup ( U_1\times \mathbb{P}_\mathbb{C}^1)$, where $(u,\xi)\in U \times \mathbb{P}_\mathbb{C}^1$ and $(u_1,\xi_1)\in U_1 \times \mathbb{P}_\mathbb{C}^1$ are identified if and only if $\xi=e^{\frac{t_0}{u_1}}\xi_1$, and $[u]=p+u_1$, where $t_0$ is a complex number not belonging to the lattice such that the ruled surface $S$ is represented by the $1$-cocycle $\eta(t_0)$. A holomorphic vector field $\theta\in H^0(S,\Theta_S)$ can be expressed in the form 
\begin{align}
\theta=(a_0(u)+a_1(u)\xi+a_2(u)\xi^2)\frac{\partial}{\partial \xi}+b(u)\frac{\partial}{\partial u}\,\,\,\,\,\text{ on $U\times \mathbb{P}_\mathbb{C}^1$},
\end{align}
 where $a_0(u),a_1(u),a_2(u)$ and $b(u)$ are holomorphic functions of $[u]\in U$. We note that $\frac{\partial}{\partial \xi_1}=e^{\frac{t_0}{u_1}}\frac{\partial}{\partial \xi}$, and $\frac{\partial}{\partial u_1}=-\frac{t_0}{u_1^2}e^{\frac{t_0}{u_1}}\xi_1\frac{\partial}{\partial \xi}+\frac{\partial}{\partial u}$. If we write $\theta$ in terms of $(u_1,\xi_1)$, we have 
\begin{align}
\theta=\left(a_0(u)e^{-\frac{t_0}{u_1}}+\left(a_1(u)+\frac{t_0}{u_1^2}b(u)\right)\xi_1+a_2(u)e^{\frac{t_0}{u_1}}\xi_1^2\right)\frac{\partial}{\partial \xi_1}+b(u)\frac{\partial}{\partial u_1}.
\end{align}
 Then $b(u):=b$ is a constant. In a neighborhood of $p$, $a_1(u)$ has the form: $a_1(u)=-\frac{t_0 b}{u_1^2}+\alpha_0+\alpha_1u_1+\cdots$, where $\alpha_i\in \mathbb{C}$ so that $a_1(u)=c-t_0b\wp(u-p)$, where $c$ is a constant and $\wp(u)$ is the Weierstrass $\wp$-function with the periods $(1,\omega)$.\footnote{For the theory of elliptic functions, we refer to \cite{Mar77} Vol. III Chapter $5$.} Since $\dim_\mathbb{C} H^0(S,\Theta_S)=2$,  we have $a_0(u)=a_2(u)=0$, and
\begin{align}\label{a11}
\xi\frac{\partial}{\partial \xi},\,\,\,\,\,-t_0\wp(u-p)\xi\frac{\partial}{\partial \xi}+\frac{\partial}{\partial u} \in H^0(S,\Theta_S)
\end{align}
forms a basis of $H^0(S,\Theta_S)$. On the other hand, $\dim_\mathbb{C} H^1(S,\Theta_S)=2$ and
\begin{align}\label{a12}
\frac{1}{u_1}\frac{\partial}{\partial u_1},\,\,\,\,\,\frac{\xi_1}{u_1}\frac{\partial}{\partial \xi_1}\in C^1(\mathcal{U}, \Theta_S)
\end{align}
forms a basis of $H^1(S,\Theta_S)$.

\subsection{Construction of $S=S_n(n\geq 1)$, and descriptions of $H^0(S,\Theta_S)$ and $H^1(S,\Theta_S)$}\label{subsection3}\

We set $S=(U\times \mathbb{P}_\mathbb{C}^1)\cup (U_1\times \mathbb{P}_\mathbb{C}^1)$, where $(u,\xi)\in U\times \mathbb{P}_\mathbb{C}^1$ and $(u_1,\xi_1)\in U_1\times \mathbb{P}_\mathbb{C}^1$ are identified if and only if $\xi=u_1^n\xi_1$ and $[u]=p+u_1$. Then $\dim_\mathbb{C}H^0(S,\Theta_S)=n+1$ and
\begin{align}\label{a16}
\xi\frac{\partial}{\partial \xi},\,\,\,\,\,\xi^2\frac{\partial}{\partial \xi},\,\,\,\,\,\wp^{(k)}(u-p)\xi^2\frac{\partial}{\partial \xi}\in H^0(S, \Theta_S),\,\,\,\,\,k=0,...,n-2,
\end{align}
forms a basis of $H^0(S,\Theta_S)$, where $\wp(u)$ is the Weierstrass $\wp$-function with the periods $(1,\omega)$, and $\wp^{(k)}(u)$ denotes the $k$-th derivative of $\wp(u)$. On the other hand, $\dim_\mathbb{C}H^1(S,\Theta_S)=n+1$ and
\begin{align}\label{a17}
\frac{1}{u_1}\frac{\partial}{\partial u_1},\,\,\,\,\,\frac{1}{u_1^{n+1}}\frac{\partial}{\partial \xi_1},\,\,\,\,\,\frac{1}{u_1^k}\frac{\partial}{\partial \xi_1}\in C^1(\mathcal{U},\Theta_S), \,\,\,\,\,k=1,...,n-1,
\end{align}
forms a basis of $H^1(S,\Theta_S)$.
\subsection{Construction of $S=A_0$, and descriptions of $H^0(S,\Theta_S)$ and $H^1(S,\Theta_S)$}\label{subsection4}\

We set $S=(U\times \mathbb{P}_\mathbb{C}^1)\cup ( U_1\times \mathbb{P}_\mathbb{C}^1)$, where $(u,\xi)\in U\times \mathbb{P}_\mathbb{C}^1$ and $(u_1,\xi_1)\in U_1\times \mathbb{P}_\mathbb{C}^1$ are identified if and only if $\xi=\xi_1+\frac{1}{u_1}$ and $[u]=p+u_1$. Then $\dim_\mathbb{C}H^0(S,\Theta_S)=2$ and
\begin{align}\label{a20}
\frac{\partial}{\partial \xi},\,\,\,\,\,-\wp(u-p)\frac{\partial}{\partial \xi}+\frac{\partial}{\partial u}\in H^0(S,\Theta_S)
\end{align}
forms a basis of $H^0(S,\Theta_S)$. On the other hand, $\dim_\mathbb{C}H^1(S,\Theta_S)=2$ and
\begin{align}\label{a21}
\frac{1}{u_1}\frac{\partial}{\partial u_1},\,\,\,\,\,-\left(\frac{\xi_1^2}{u_1}+\frac{\xi_1}{u_1^2} \right)\frac{\partial}{\partial \xi_1}\in C^1(\mathcal{U},\Theta_S)
\end{align}
forms a basis of $H^1(S,\Theta_S)$.

\subsection{Construction of $S=A_{-1}$, and descriptions of $H^0(S,\Theta_S)$ and $H^1(S,\Theta_S)$}\label{subsection5}\

We set $S=(U\times \mathbb{P}_\mathbb{C}^1) \cup (U_1\times \mathbb{P}_\mathbb{C}^1)$, where $(u,\xi)\in U\times \mathbb{P}_\mathbb{C}^1$ and $(u_1,\xi_1)\in U_1\times \mathbb{P}_\mathbb{C}^1$ are identified if and only if $\xi=u_1\xi_1+\frac{1}{u_1}$ and $[u]=p+u_1$. Then $\dim_\mathbb{C}H^0(S,\Theta_S)=1$ and
\begin{align}
-3\wp(u-p)\frac{\partial}{\partial \xi}+\xi^2\frac{\partial}{\partial \xi}+2\frac{\partial}{\partial u}\in H^0(S,\Theta_S)
\end{align}
forms a basis of $H^0(S,\Theta_S)$. On the other hand, $\dim_\mathbb{C}H^1(S,\Theta_S)=1$ and
\begin{align}
\frac{1}{u_1}\frac{\partial}{\partial u_1}\in C^1(\mathcal{U},\Theta_S)
\end{align}
forms a basis of $H^1(S,\Theta_S)$.

\section{Holomorphic Poisson structures on ruled surfaces $S$ over an elliptic curve, and descriptions of $H^1(S,\wedge^2 \Theta_S)$}\label{section2}

In this section, based on the biholomorphic classification of ruled surfaces over an elliptic curve in the previous section, we explicitly describe holomorphic Poisson structures on ruled surfaces over an elliptic curve. We remark that the dimensions of holomorphic Poisson structures $H^0(S,\wedge^2 \Theta_S)$ on ruled surfaces $S$ over an elliptic curve were computed in \cite{BM05}. We note that as the Chern numbers $c_1^2$ and $c_2$ vanish, and $H^2(S,\wedge^2 \Theta_S)=0$, by  Hirzebruch-Riemann-Roch theorem, we have $\dim_\mathbb{C} H^1(S,\wedge^2 \Theta_S)=\dim_\mathbb{C} H^0(S,\wedge^2 \Theta_S)$.

We keep the notations in section $\ref{section1}$. Then a bivector field on $U\times \mathbb{P}_\mathbb{C}^1$ is of the form
\begin{align}\label{a2}
\left( a_0(u)+a_1(u)\xi+a_2(u)\xi^2\right)\frac{\partial}{\partial \xi}\wedge \frac{\partial}{\partial u},\,\,\,\,\,\text{on $\,\,\,U\times \mathbb{P}_\mathbb{C}^1$},
\end{align}
where $a_0(u),a_1(u)$ and $a_2(u)$ are holomorphic functions of $[u]\in U$. On the other hand, a bivector field on $U_1\times \mathbb{P}_\mathbb{C}^1$ is of the form
\begin{align}\label{a3}
\left( a_{10}(u_1)+a_{11}(u_1)\xi_1+a_{12}(u_1)\xi_1^2\right)\frac{\partial}{\partial \xi_1}\wedge \frac{\partial}{\partial u_1},\,\,\,\,\,\text{on $\,\,\,U_1\times \mathbb{P}_\mathbb{C}^1$},
\end{align}
where $a_{10}(u_1),a_{11}(u_1)$ and $a_{12}(u_1)$ are holomorphic functions of $u_1\in U_1$.
\subsection{Holomorphic Poisson structures on $S=S_0=X\times \mathbb{P}_\mathbb{C}^1$, and descriptions of $H^1(S,\wedge^2 \Theta_S)$}\

We keep the notations in subsection $\ref{subsection1}$. Then $\dim_\mathbb{C}H^0(S,\wedge^2\Theta_S)=3$ and
\begin{align}\label{a4}
\frac{\partial}{\partial \xi}\wedge \frac{\partial}{\partial u},\,\,\,\,\,\xi\frac{\partial}{\partial \xi}\wedge \frac{\partial}{\partial u},\,\,\,\,\,\xi^2\frac{\partial}{\partial \xi}\wedge \frac{\partial}{\partial u}\in H^0(S,\wedge^2 \Theta_S)
\end{align}
forms a basis of $H^0(S,\wedge^2 \Theta_S)$. On the other hand, since $\dim_\mathbb{C}H^1(S,\wedge^2\Theta_S)=3$ and there is no elliptic function of order $1$, 
\begin{align}\label{a5}
\frac{1}{u_1}\frac{\partial}{\partial \xi_1}\wedge \frac{\partial}{\partial u_1},\,\,\,\,\,\frac{\xi_1}{u_1}\frac{\partial}{\partial \xi_1}\wedge \frac{\partial}{\partial u_1},\,\,\,\,\,\frac{\xi^2}{u_1}\frac{\partial}{\partial \xi_1}\wedge \frac{\partial}{\partial u_1}\in C^1(\mathcal{U},\wedge^2 \Theta_S)
\end{align}
forms a basis of $H^1(S, \wedge^2\Theta_S)$.
\subsection{Holomorphic Poisson structures on $S\in \mathcal{S}_0,S\ne S_0$, and descriptions of $H^1(S,\wedge^2 \Theta_S)$}\

We keep the notations in subsection $\ref{subsection2}$. Then $\frac{\partial}{\partial \xi_1}\wedge \frac{\partial}{\partial u_1}=e^{\frac{t_0}{u_1}}\frac{\partial}{\partial \xi}\wedge \frac{\partial}{\partial u}$. If we write a holomorphic bivector field $\Lambda_0\in H^0(S,\wedge^2\Theta_S)$ of the form $(\ref{a2})$ in terms of $(u_1,\xi_1)$, we have
\begin{align*}
\Lambda_0=\left(a_0(u)e^{-\frac{t_0}{u_1}}+a_1(u)\xi_1+a_2(u)e^{\frac{t_0}{u_1}}\xi_1^2 \right)\frac{\partial}{\partial \xi_1}\wedge \frac{\partial}{\partial u_1}
\end{align*}
Then $a_0(u)=a_2(u)=0$ and $a_1(u):=a_1$ is a constant. Hence $\dim_\mathbb{C}H^0(S,\wedge^2\Theta_S)=1$ and
\begin{align}\label{a9}
\xi\frac{\partial}{\partial \xi}\wedge \frac{\partial}{\partial u}\in H^0(S,\wedge^2 \Theta_S)
\end{align}
forms a basis of $H^0(S,\wedge^2 \Theta_S)$. On the other hand, since $\dim_\mathbb{C}H^1(S,\wedge^2 \Theta_S)=1$, and there is no elliptic function of order $1$, we see that
\begin{align}\label{a10}
\frac{\xi_1}{u_1}\frac{\partial}{\partial \xi_1}\wedge \frac{\partial}{\partial u_1}\in C^1(\mathcal{U},\wedge^2 \Theta_S)
\end{align}
forms a basis of $H^1(S,\wedge^2 \Theta_S)$.
\subsection{Holomorphic Poisson structures on $S=S_n(n\geq 1)$, and descriptions of $H^1(S,\wedge^2 \Theta_S)$}\

We keep the notations in subsection $\ref{subsection3}$. Then $\frac{\partial}{\partial \xi_1}=u_1^n\frac{\partial}{\partial \xi}$ and $\frac{\partial}{\partial u_1}=nu_1^{n-1}\xi_1\frac{\partial}{\partial u}+\frac{1}{u}$ and so $\frac{\partial}{\partial \xi_1}\wedge \frac{\partial}{\partial u_1}=u_1^n \frac{\partial}{\partial \xi}\wedge \frac{\partial}{\partial u}$. If we write a holomorphic bivector field $\Lambda_0\in H^0(S,\wedge^2\Theta_S)$ of the form $(\ref{a2})$ in terms of $(u_1,\xi_1)$, we have
\begin{align*}
\Lambda_0=\left( a_0(u)\frac{1}{u_1^n}+a_1(u)\xi_1+a_2(u)u_1^n \xi_1^2\right)\frac{\partial}{\partial \xi_1}\wedge \frac{\partial}{\partial u_1}.
\end{align*}
It follows that $a_0(u)=0$ and $a_1(u):=a_1$ is a constant. Since $p$ is a pole of $a_2(u)$ of order at most $n$, 
\begin{align*}
a_2(u)=c_0+c_1\wp(u-p)+c_2\wp'(u-p)+\cdots +c_{n-1}\wp^{(n-2)}(u-p), \,\,\,\,\,c_0,c_1,...,c_{n-1}\in \mathbb{C},
\end{align*}
where $\wp(u)$ is the Weierstrass $\wp$-function with the periods $(1,\omega)$, and $\wp^{(k)}(u)$ denotes the $k$-th derivative of $\wp(u)$. Hence $\dim_\mathbb{C}H^0(S,\wedge^2\Theta_S)=n+1$ and
\begin{align}\label{a13}
\xi\frac{\partial}{\partial \xi}\wedge \frac{\partial}{\partial u},\,\,\,\,\,\xi^2 \frac{\partial}{\partial \xi}\wedge \frac{\partial}{\partial u},\,\,\,\,\,\wp^{(k)}(u-p)\xi^2\frac{\partial}{\partial \xi}\wedge \frac{\partial}{\partial u}\in H^0(S,\wedge^2 \Theta_S),\,\,\,\,\,k=0,...,n-2,
\end{align}
forms a basis of $H^0(S,\wedge^2 \Theta_S)$. On the other hand, we note that $\dim_\mathbb{C}H^1(S,\wedge^2 \Theta_S)=n+1$. We claim that
\begin{align}\label{a14}
\frac{\xi_1}{u_1}\frac{\partial}{\partial \xi_1}\wedge\frac{\partial}{\partial u_1},\,\,\,\,\,\frac{1}{u_1^{n+1}}\frac{\partial}{\partial \xi_1}\wedge \frac{\partial}{\partial u_1},\,\,\,\,\,\frac{1}{u_1^k}\frac{\partial}{\partial \xi_1}\wedge \frac{\partial}{\partial u_1}\in C^1(\mathcal{U},\wedge^2\Theta_S),\,\,\,\,\,k=1,...,n-1,
\end{align}
forms a basis of $H^1(S,\wedge^2 \Theta_S)$. Indeed, assume that $\frac{d}{u_1}\xi_1\frac{\partial}{\partial \xi_1}\wedge\frac{\partial}{\partial u_1}+\left(\frac{c_1}{u_1^{n+1}}+ \frac{c_2}{u_1}+\frac{c_3}{u_1^2}+\cdots+\frac{c_n}{u_1^{n-1}}\right)\frac{\partial}{\partial \xi_1}\wedge \frac{\partial}{\partial u_1}$ is written as a difference of two bivector fields of the forms $(\ref{a2})$ and $(\ref{a3})$ for some constants $d,c_1,c_2,...,c_{n-1}$:
\begin{align*}
&\frac{d}{u_1}\xi_1\frac{\partial}{\partial \xi_1}\wedge\frac{\partial}{\partial u_1}+\left(\frac{c_1}{u_1^{n+1}}+ \frac{c_2}{u_1}+\frac{c_3}{u_1^2}+\cdots+\frac{c_n}{u_1^{n-1}}\right)\frac{\partial}{\partial \xi_1}\wedge \frac{\partial}{\partial u_1}\\
&=\left( a_{10}(u_1)-a_0(u)\frac{1}{u_1^n}+(a_{10}(u_1)-a_1(u))\xi_1+(a_{20}(u_1)-a_2(u)u_1^n) \xi_1^2\right)\frac{\partial}{\partial \xi_1}\wedge \frac{\partial}{\partial u_1}.
\end{align*}
This implies that $a_1(u)=-\frac{d}{u_1}+a_{10}(u_1)$. Since there is no elliptic function of order $1$, we have $d=0$. On the other hand, $\frac{a_0(u)}{u_1^n}=-\left(\frac{c_1}{u_1^{n+1}}+ \frac{c_2}{u_1}+\frac{c_3}{u_1^2}+\cdots+\frac{c_n}{u_1^{n-1}}\right)+a_{10}(u_1)$ and hence in a neighborhood of $u_1=0$, we have
$a_0(u)=-\frac{c_1}{u_1}-c_nu_1-\cdots -c_2u_1^{n-1}+\alpha_0u_1^n+\alpha_1u_1^{n+1}+\cdots$, where $\alpha_i\in \mathbb{C}$. Hence $c_1=0$ and $a_0(u)$ is a constant so that $a_0(u)=0$ and $c_2=\cdots=c_n=0$. This proves the claim.

\subsection{Holomorphic Poisson structures on $S=A_0$, and descriptions of $H^1(S,\wedge^2 \Theta_S)$}\

We keep the notations in subsection $\ref{subsection4}$. Then $\frac{\partial}{\partial \xi_1}=\frac{\partial}{\partial \xi}$ and $\frac{\partial}{\partial u_1}=-\frac{1}{u_1^2}\frac{\partial}{\partial \xi}+\frac{\partial}{\partial u}$ and $\frac{\partial}{\partial \xi_1}\wedge \frac{\partial}{\partial u_1}=\frac{\partial}{\partial \xi}\wedge \frac{\partial}{\partial u}$. If we write a holomorphic bivector field $\Lambda_0\in H^0(S,\wedge^2\Theta_S)$ of the form $(\ref{a2})$ in terms of $(u_1,\xi_1)$, we have
\begin{align*}
\Lambda_0
=\left( a_0(u)+\frac{a_1(u)}{u_1}+\frac{a_2(u)}{u_1^2}+\left(a_1(u)+2a_2(u)\frac{1}{u_1}\right) \xi_1+ a_2(u)\xi_1^2\right)\frac{\partial}{\partial \xi_1}\wedge \frac{\partial}{\partial u_1}
\end{align*}
Then $a_2(u):=a_2$ is a constant. In a neighborhood of $p$, $a_1(u)$ has the following form:
\begin{align*}
a_1(u)=-\frac{2a_2}{u_1}+\alpha_0+\alpha_1u_1+\alpha_2u_1^2+\cdots,\,\,\,\,\,\alpha_1\in \mathbb{C}.
\end{align*}
As there exists no elliptic function of order $1$, we have $a_2=0$ and $a_1(u):=a_1$ is a constant. In a neighborhood of $p$, $a_0(u)$ has the following form:
\begin{align*}
a_0(u)=-\frac{a_1}{u_1}+\beta_0+\beta_1u_1+\beta_2u_1^2+\cdots, \,\,\,\,\,\beta_i\in \mathbb{C}.
\end{align*}
Since there exists no elliptic function of order $1$, we have $a_1=0$ and $a_0(u):=a_0$ is a constant: $a_0(u)=a_0$. Hence $\dim_\mathbb{C}H^0(S,\wedge^2\Theta_S)=1$ and
\begin{align}\label{a18}
\frac{\partial}{\partial \xi}\wedge \frac{\partial}{\partial u}\in H^0(S,\wedge^2 \Theta_S)
\end{align}
forms a basis of $H^0(S,\wedge^2 \Theta_S)$. On the other hand, we note that $\dim_\mathbb{C} H^1(S,\wedge^2 \Theta_S)=1$, and there is no elliptic function of order $1$, we see that 
\begin{align}\label{a19}
\frac{\xi_1^2}{u_1}\frac{\partial}{\partial \xi_1}\wedge \frac{\partial}{\partial u_1}\in C^1(\mathcal{U},\wedge^2 \Theta_S)
\end{align}
forms a basis of $H^1(S,\wedge^2 \Theta_S)$.

\subsection{Holomorphic Poisson structures on $S=A_{-1}$, and descriptions of $H^1(S,\wedge^2 \Theta_S)$}\

We keep the notations in subsection $\ref{subsection5}$. Then $\frac{\partial}{\partial \xi_1}=u_1\frac{\partial}{\partial \xi},\frac{\partial}{\partial u_1}=\left(\xi_1-\frac{1}{u_1^2}\right)\frac{\partial}{\partial \xi}+\frac{\partial}{\partial u}$ and $\frac{\partial}{\partial \xi_1}\wedge \frac{\partial}{\partial u_1}= u_1\frac{\partial}{\partial \xi_1}\wedge \frac{\partial}{\partial u_1}$. If we write a holomorphic bivector field $\Lambda_0\in H^0(S,\wedge^2\Theta_S)$ of the form $(\ref{a2})$ in terms of $(u_1,\xi_1)$, we have
\begin{align*}
\Lambda_0&=\left(\frac{a_0(u)}{u_1}+\frac{a_1(u)}{u_1^2}+\frac{a_2(u)}{u_1^3}+\left(a_1(u)+\frac{2a_2(u)}{u_1} \right) \xi_1 +  a_2(u)u_1\xi_1^2 \right)\frac{\partial}{\partial \xi_1}\wedge \frac{\partial}{\partial u_1}.
\end{align*}
Then $a_2(u):=a_2$ is a constant. In a neighborhood of $p$, $a_1(u)$ has the following form
\begin{align*}
a_1(u)=-\frac{2a_2}{u_1}+\alpha_0+\alpha_1u_1+\alpha_2u_1^2+\cdots,\,\,\,\,\,\alpha_i\in \mathbb{C}.
\end{align*}
Hence $a_2=0$ and $a_1(u):=a_1$ is a constant. In a neighborhood of $p$, $a_0(u)$ has the following form:
\begin{align*}
a_0(u)=-\frac{a_1}{u_1}+\beta_0u_1+\beta_1u_1^2+\cdots, \,\,\,\,\,\,\beta_i\in \mathbb{C}.
\end{align*}
Hence $a_1=0$ and $a_0 (u) =0$ so that we have
\begin{align}\label{a22}
 \text{$H^0(S,\wedge^2 \Theta_S)=0\,\,\,\,$ and $\,\,H^1(S,\wedge^2 \Theta_S)=0$}.
\end{align}

\section{Computations of $\mathbb{H}^0(S, \Theta_S^\bullet), \mathbb{H}^1(S,\Theta_S^\bullet)$ and $\mathbb{H}^2(S,\Theta_S^\bullet)$ for Poisson ruled surfaces $(S,\Lambda_0)$ over an elliptic curve}\label{section3}

In this section, given a Poisson ruled surface $(S,\Lambda_0)$ over an elliptic curve $X$ described in section $\ref{section2}$, we compute cohomology groups $\mathbb{H}^0(S,\Theta_S^\bullet), \mathbb{H}^1(S,\Theta_S^\bullet)$ and $\mathbb{H}^2(S,\Theta_S^\bullet)$ of the complex of sheaves $\Theta_S^\bullet:\Theta_S \xrightarrow{[\Lambda_0,-]} \wedge^2\Theta_S\xrightarrow{[\Lambda_0,-]} \wedge^3 \Theta_S=0$. By considering the spectral sequence associated with $\Theta_S\xrightarrow{[\Lambda_0,-]} \wedge^2 \Theta_S\to 0$ and $H^2(S,\Theta_S)=0$, we have 

\begin{lemma}\label{a25}
Let $(S,\Lambda_0)$ be a Poisson ruled surface over an elliptic curve $X$. Then 
\begin{align*}
\mathbb{H}^0(S, \Theta_{S}^\bullet)&\cong ker(H^0(S, \Theta_{S})\xrightarrow{[\Lambda_0,-]} H^0(S, \wedge^2 \Theta_{S}))\\
\mathbb{H}^1(S, \Theta_{S}^\bullet)&\cong coker(H^0(S, \Theta_{S})\xrightarrow{[\Lambda_0,-]} H^0(S, \wedge^2 \Theta_{S}))\oplus ker(H^1(S, \Theta_{S})\xrightarrow{[\Lambda_0,-]} H^1(S,\wedge^2 \Theta_{S}))\\
\mathbb{H}^2(S,\Theta_{S}^\bullet)&\cong coker(H^1(S,\Theta_{S})\xrightarrow{[\Lambda_0,-]}H^1(S,\wedge^2 \Theta_{S}))
\end{align*}
\end{lemma}

In the following, we keep the notations in section $\ref{section1}$ and section $\ref{section2}$.

\subsection{The case of $S=S_0=X\times \mathbb{P}_\mathbb{C}^1$}\

From $(\ref{a4})$, let $(S,\Lambda_0)$ be a Poisson ruled surface over an elliptic curve $X$ given by
\begin{align}\label{a6}
\Lambda_0=(A+B\xi+C\xi^2)\frac{\partial}{\partial \xi}\wedge \frac{\partial}{\partial u}=(A+B\xi_1+C\xi_1^2)\frac{\partial}{\partial \xi_1}\wedge \frac{\partial}{\partial u_1},
\end{align}
where $(A,B,C)\in \mathbb{C}^3$ are constants. Let us compute $H^0(S, \Theta_{S})\xrightarrow{[\Lambda_0,-]} H^0(S, \wedge^2 \Theta_{S})$. From $(\ref{a7})$ and $(\ref{a4})$, we have
\begin{align*}
&[(A+B\xi+C\xi^2)\frac{\partial}{\partial \xi}\wedge \frac{\partial}{\partial u}, t_0\frac{\partial}{\partial u}+(t_1+t_2\xi+t_3\xi^2)\frac{\partial}{\partial \xi}]\\
&=\left((-Bt_1+At_2)+2(-Ct_1+At_3)\xi+(-Ct_2+Bt_3)\xi^2 \right)\frac{\partial}{\partial \xi}\wedge \frac{\partial}{\partial u}
\end{align*}
where $t_0,t_1,t_2,t_3\in \mathbb{C}$. If $(A,B,C)=0$, then $\dim_\mathbb{C}\mathbb{H}^0(S,\Theta_S^\bullet)=4$, and if $(A,B,C)\ne 0$, then $\dim_\mathbb{C} \mathbb{H}^0(S,\Theta_S^\bullet)=2$ by Lemma $\ref{a25}$. Let us compute $H^1(S,\Theta_{S})\xrightarrow{[\Lambda_0,-]}H^1(S,\wedge^2 \Theta_{S})$. From $(\ref{a8})$ and $(\ref{a5})$, we have

\begin{align*}
&[(A+B\xi_1+C\xi_1^2)\frac{\partial}{\partial \xi_1}\wedge \frac{\partial}{\partial u_1}, \frac{t_0}{u_1}\frac{\partial}{\partial u_1}+\frac{t_1+t_2\xi_1+t_3\xi_1^2}{u_1}\frac{\partial}{\partial \xi_1}]\\
&=\left((-Bt_1+At_2)\frac{1}{u_1}+2(-Ct_1+At_3)\frac{\xi_1}{u_1}+(-Ct_2+Bt_3)\frac{\xi_1^2}{u_1}-(A+B\xi_1+C\xi_1^2)\frac{t_0}{u_1^2} \right)\frac{\partial}{\partial \xi_1}\wedge \frac{\partial}{\partial u_1}\\
&\equiv \left((-Bt_1+At_2)\frac{1}{u_1}+2(-Ct_1+At_3)\frac{\xi_1}{u_1}+(-Ct_2+Bt_3)\frac{\xi_1^2}{u_1}\right)\frac{\partial}{\partial \xi_1}\wedge \frac{\partial}{\partial u_1}
\end{align*}
where $t_0,t_1,t_2,t_3\in \mathbb{C}$. Here $a\equiv b$ means that $a$ and $b$ represent the same cohomology class in $H^1(S,\wedge^2 \Theta_S)$. If $(A,B,C)=0$, then $\dim_\mathbb{C} \mathbb{H}^1(S,\Theta_S^\bullet)=7$, and $\dim_\mathbb{C}\mathbb{H}^2(S,\Theta_S^\bullet)=3$ by Lemma $\ref{a25}$. If $(A,B,C)\ne 0$, then $\dim_\mathbb{C}\mathbb{H}^1(S,\Theta_S^\bullet)=3$ and $\dim_\mathbb{C}\mathbb{H}^2(S,\Theta_S^\bullet)=1$ by Lemma $\ref{a25}$. In this case, choose $(F_1,F_2,F_2)\ne 0\in \mathbb{C}^3$ such that $(F_1+F_2\xi+F_3\xi^2)\frac{\partial}{\partial \xi}\wedge \frac{\partial}{\partial u}$ is not in the image of $H^0(S, \Theta_S )\xrightarrow{[\Lambda_0,-]} H^0(S,\wedge^2 \Theta_S )$. Then we have
\begin{align}\label{a35}
\mathbb{H}^1(S,\Theta_S^\bullet)\cong \left\langle (F_1+F_2\xi+F_3\xi^2)\frac{\partial}{\partial \xi}\wedge \frac{\partial}{\partial u} \right\rangle \oplus \left\langle \frac{1}{u_1}\frac{\partial}{\partial u_1},(A+B\xi_1+C\xi_1^2)\frac{1}{u_1}\frac{\partial}{\partial \xi_1} \right\rangle
\end{align}
We summarize our computation in Table $\ref{ruled}$.

\subsection{The case of $S\in \mathcal{S}_0, S\ne S_0$}\

From $(\ref{a9})$, let $(S,\Lambda_0)$ be a Poisson ruled surface over an elliptic curve $X$ given by
\begin{align}\label{a36}
\Lambda_0=A\xi\frac{\partial}{\partial \xi}\wedge \frac{\partial}{\partial u}=A\xi_1\frac{\partial}{\partial \xi_1}\wedge \frac{\partial}{\partial u_1},
\end{align}
where $A$ is a constant. Let us compute $H^0(S, \Theta_{S})\xrightarrow{[\Lambda_0,-]} H^0(S, \wedge^2 \Theta_{S})$. From $(\ref{a11})$ and $(\ref{a9})$, we have
\begin{align*}
[A\xi\frac{\partial}{\partial \xi}\wedge \frac{\partial}{\partial u},c_0\frac{\partial}{\partial u}+\left(c_1-c_0t_0\wp(u-p)\right)\xi\frac{\partial}{\partial \xi}]=0,\,\,\,\,\,c_0,c_1\in \mathbb{C}.
\end{align*}
Hence $\dim_\mathbb{C}\mathbb{H}^0(S,\Theta_S^\bullet)=2$ by Lemma $\ref{a25}$. Let us compute $H^1(S, \Theta_{S})\xrightarrow{[\Lambda_0,-]} H^1(S, \wedge^2 \Theta_{S})$. From $(\ref{a12})$ and $(\ref{a10})$,
\begin{align*}
[A\xi_1\frac{\partial}{\partial \xi_1}\wedge \frac{\partial}{\partial u_1}, c_0\frac{1}{u_1}\frac{\partial}{\partial u_1}+c_1\frac{\xi_1}{u_1}\frac{\partial}{\partial \xi_1}]=-Ac_0\frac{\xi_1}{u_1^2}\frac{\partial}{\partial \xi_1}\wedge \frac{\partial}{\partial u_1}\equiv 0, \,\,\,\,\,c_0,c_1\in \mathbb{C}.
\end{align*}
Here $a\equiv b$ means that $a$ and $b$ represent the same cohomology class in $H^1(S,\wedge^2 \Theta_S)$. Hence $\dim_\mathbb{C} \mathbb{H}^1(S,\Theta_S^\bullet)=3$, $\dim_\mathbb{C}\mathbb{H}^2(S,\Theta_S^\bullet)=1$ by Lemma $\ref{a25}$, and we have
\begin{align}\label{a37}
\mathbb{H}^1(S,\Theta_S^\bullet)\cong \left\langle \xi\frac{\partial}{\partial \xi}\wedge \frac{\partial}{\partial u} \right \rangle \oplus \left\langle \frac{1}{u_1}\frac{\partial}{\partial u_1}, \frac{\xi_1}{u_1}\frac{\partial}{\partial \xi_1}  \right\rangle
\end{align}
We summarize our computations in Table $\ref{ruled}$.
\subsection{The case of $S=S_n(n\geq 1)$}\

From $(\ref{a13})$, let $(S,\Lambda_0)$ be a Poisson ruled surface over an elliptic curve $X$ given by
\begin{align}\label{a15}
\Lambda_0=(a_0\xi+A(u)\xi^2)\frac{\partial}{\partial \xi}\wedge \frac{\partial}{\partial u}
         =(a_0\xi_1+A(u)u_1^n \xi_1^2)\frac{\partial}{\partial \xi_1}\wedge \frac{\partial}{\partial u_1}.
\end{align}
where $A(u):=c_0+c_1\wp(u-p)+c_2\wp'(u-p)+\cdots +c_{n-1}\wp^{(n-2)}(u-p)$, and $a_0,c_0,c_1,\cdots,c_{n-1}$ are constants. Let us compute $H^0(S, \Theta_{S})\xrightarrow{[\Lambda_0,-]} H^0(S, \wedge^2 \Theta_{S})$. From $(\ref{a16})$ and $(\ref{a13})$, setting $B(u):=t_0+t_1\wp(u-p)+c_2\wp'(u-p)+\cdots +t_{n-1}\wp^{(n-2)}(u-p)$,
\begin{align*}
&[(a_0\xi+A(u)\xi^2)\frac{\partial}{\partial \xi}\wedge \frac{\partial}{\partial u}, (b_0\xi+B(u)\xi^2)\frac{\partial}{\partial \xi}]
=(a_0B(u)-b_0A(u)) \xi^2\frac{\partial}{\partial \xi}\wedge \frac{\partial}{\partial u},\,\,\,\,\,b_0,t_0,t_1,...,t_{n-1}\in \mathbb{C}.
\end{align*}
Then $ker(H^0(S, \Theta_{S})\xrightarrow{[\Lambda_0,-]} H^0(S, \wedge^2 \Theta_{S}) )$ is given by $a_0B(u)=b_0A(u)$ if and only if
\begin{align*}
a_0t_0&=b_0c_0, \\
a_0t_1&=b_0c_1\\
&\cdots\\
 a_0t_{n-1}&=b_0c_{n-1}
\end{align*}

 If $a_0\ne 0$, then $t_0=\frac{b_0c_0}{a_0}, t_1=\frac{b_0c_1}{a_0},..., t_{n-1}=\frac{b_0 c_{n-1}}{a_0}$. Hence $\dim_\mathbb{C} \mathbb{H}^0(S,\Theta_S^\bullet)=1$ and $\dim_\mathbb{C}coker( H^0(S, \Theta_{S})\xrightarrow{[\Lambda_0,-]} H^0(S, \wedge^2 \Theta_{S}))=1$ by Lemma $\ref{a25}$. In this case,
\begin{align}\label{a30}
coker(H^0(S, \Theta_{S})\xrightarrow{[\Lambda_0,-]} H^0(S, \wedge^2 \Theta_{S}) )\cong \left\langle\xi\frac{\partial}{\partial \xi}\wedge \frac{\partial}{\partial u} \right\rangle. 
\end{align}

 If $a_0= 0$, and $A(u)=0$, then $\dim_\mathbb{C} \mathbb{H}^0(S,\Theta_S^\bullet)=n+1$, and 
 \begin{align}\label{b1}
 \dim_\mathbb{C}coker(H^0(S, \Theta_{S})\xrightarrow{[\Lambda_0,-]} H^0(S, \wedge^2 \Theta_{S}))=n+1.
 \end{align}

 If $a_0=0$, and $A(u)\ne 0$, then $ker(H^0(S, \Theta_{S})\xrightarrow{[\Lambda_0,-]} H^0(S, \wedge^2 \Theta_{S}))$ is given by $-b_0A(u)=0$ if and only if $b_0=0$. Hence $\dim_\mathbb{C} \mathbb{H}^0(S,\Theta_S^\bullet)= n$, and 
\begin{align}\label{b2}
 \dim_\mathbb{C}coker(H^0(S, \Theta_{S})\xrightarrow{[\Lambda_0,-]} H^0(S, \wedge^2 \Theta_{S}))=n
\end{align} 

Next let us compute $H^1(S, \Theta_{S})\xrightarrow{[\Lambda_0,-]} H^1(S, \wedge^2 \Theta_{S})$. From $(\ref{a17})$ and $(\ref{a14})$,
\begin{align*}
&[(a_0\xi_1+A(u)u_1^n \xi_1^2)\frac{\partial}{\partial \xi_1}\wedge \frac{\partial}{\partial  u_1},\frac{b}{u_1}\frac{\partial}{\partial u_1}+\left( \frac{t_1}{u_1^{n+1}}+\frac{t_2}{u_1}+\frac{t_3}{u_1^2}+\cdots +\frac{t_n}{u_1^{n-1}}\right)\frac{\partial}{\partial \xi_1}]\\
&=\left(-\frac{b}{u_1}\frac{\partial (A(u)u_1^n )}{\partial u_1}\xi_1^2-\left( \frac{t_1}{u_1^{n+1}}+\frac{t_2}{u_1}+\frac{t_3}{u_1^2}+\cdots +\frac{t_n}{u_1^{n-1}}\right)(a_0+2A(u)u_1^n\xi_1)-(a_0\xi+A(u)u_1^n\xi_1^2)\frac{b}{u_1^2}\right)\frac{\partial}{\partial \xi_1}\wedge \frac{\partial}{\partial u_1}\\
&\equiv -a_0\left( \frac{t_1}{u_1^{n+1}}+\frac{t_2}{u_1}+\frac{t_3}{u_1^2}+\cdots +\frac{t_n}{u_1^{n-1}}\right)\frac{\partial}{\partial \xi_1}\wedge \frac{\partial}{\partial u_1} \\
&-2\left( \frac{t_1}{u_1}+t_2u_1^{n-1}+{t_3}u_1^{n-2}+\cdots +t_nu_1\right) \left(c_0+c_1\wp(u-p)+c_2\wp'(u-p)+\cdots +c_{n-1}\wp^{(n-2)}(u-p)\right)\xi_1\frac{\partial}{\partial \xi_1}\wedge \frac{\partial}{\partial u_1}\\
&\equiv -a_0\left( \frac{t_1}{u_1^{n+1}}+\frac{t_2}{u_1}+\frac{t_3}{u_1^2}+\cdots +\frac{t_n}{u_1^{n-1}}\right)\frac{\partial}{\partial \xi_1}\wedge \frac{\partial}{\partial u_1} -2(t_1c_0+ t_2c_{n-1}+t_3c_{n-2}+\cdots+t_nc_1)\frac{1}{u_1}\xi_1\frac{\partial}{\partial \xi_1}\wedge \frac{\partial}{\partial u_1}
\end{align*}
where $b,t_1,t_2,...,t_n\in \mathbb{C}$. Here $a\equiv b$ means that $a$ and $b$ represent the same cohomology class in $H^1(S,\wedge^2 \Theta_S)$.

If $a_0\ne 0$, then $ker(H^1(S, \Theta_{S})\xrightarrow{[\Lambda_0,-]} H^1(S, \wedge^2 \Theta_{S}))$ is given by $t_1=\cdots=t_n=0$. Hence $\mathbb{H}^2(S,\Theta_S^\bullet)=n+1-n=1$ by Lemma $\ref{a25}$. From $(\ref{a30})$, $\dim_\mathbb{C} \mathbb{H}^1(S,\Theta_S^\bullet)=2$ and
\begin{align}\label{a39}
\mathbb{H}^1(S,\Theta_S^\bullet)\cong \left\langle\xi\frac{\partial}{\partial \xi}\wedge \frac{\partial}{\partial u} \right\rangle \oplus \left\langle \frac{1}{u_1}\frac{\partial}{\partial u_1} \right\rangle. 
\end{align} 

If $a_0=0$ and $A(u)=0$, then $\dim_\mathbb{C}\mathbb{H}^2(S,\Theta_S^\bullet)=n+1$, and from $(\ref{b1})$,  $\dim_\mathbb{C} \mathbb{H}^1(S,\Theta_S^\bullet)=2n+2$ by Lemma $\ref{a25}$.

If $a_0=0$ and $ A(u)\ne 0$, then $\dim_\mathbb{C} \mathbb{H}^2(S,\Theta_S^\bullet)=n$, and from $(\ref{b2})$, $\dim_\mathbb{C}\mathbb{H}^1(S,\Theta_S^\bullet)=2n$ by Lemma $\ref{a25}$. In particular, if $a_0=0, n=1$ (i.e. $S=S_1$) and $A(u)=c_0\ne 0$, then
\begin{align}\label{a90}
\mathbb{H}^1(S_1,\Theta_{S_1}^\bullet)\cong \left\langle\xi\frac{\partial}{\partial \xi}\wedge \frac{\partial}{\partial u} \right\rangle \oplus \left\langle \frac{1}{u_1}\frac{\partial}{\partial u_1} \right\rangle. 
\end{align}

For later use, we note that if $\Lambda_0=0$ (i.e $a_0=0,A(u)=0$) for $n \geq 1$ or $a_0=0,A(u)\ne 0$ for $n\geq 2$, then there exist constants $(t_1,...t_n)\ne 0\in \mathbb{C}^n$ such that 
\begin{align}\label{a51}
\frac{t_1}{u_1^{n+1}}+\frac{t_2}{u_1}+\frac{t_3}{u_1^2}+\cdots +\frac{t_n}{u_1^{n-1}}\in C^1(\mathcal{U},\Theta_S)
\end{align}
is in $ker(H^1(S,\Theta_S)\xrightarrow{[\Lambda_0,-]} H^1(S,\wedge^2 \Theta_S))$.

We summarize our computation in Table $\ref{ruled}$.

\subsection{The case of $S=A_0$}\

From $(\ref{a18})$, let $(S,\Lambda_0)$ be a Poisson ruled surface over an elliptic curve $X$ given by
\begin{align}\label{a53}
\Lambda_0=a_0\frac{\partial}{\partial \xi}\wedge \frac{\partial}{\partial u}=a_0\frac{\partial}{\partial \xi_1}\wedge \frac{\partial}{\partial u_1}
\end{align}
where $a_0$ is a constant. Let us compute $ker(H^0(S, \Theta_{S})\xrightarrow{[\Lambda_0,-]} H^0(S, \wedge^2 \Theta_{S}))$. From $(\ref{a20})$ and $(\ref{a18})$,
\begin{align*}
[a_0\frac{\partial}{\partial \xi}\wedge \frac{\partial}{\partial u}, t_2\frac{\partial}{\partial u}+( t_1-t_2\wp(u-p))\frac{\partial}{\partial \xi}]=0,\,\,\,\,\,t_1,t_2\in \mathbb{C}.
\end{align*}
Hence $\dim_\mathbb{C} \mathbb{H}^0(S,\Theta_S^\bullet)=2$ by Lemma $\ref{a25}$. On the other hand, let us compute $H^1(S, \Theta_{S})\xrightarrow{[\Lambda_0,-]} H^1(S, \wedge^2 \Theta_{S})$. From $(\ref{a21})$ and $(\ref{a19})$,
\begin{align*}
[a_0\frac{\partial}{\partial \xi_1}\wedge \frac{\partial}{\partial u_1}, \frac{t_1}{u_1}\frac{\partial}{\partial u_1} -t_2\left(\frac{\xi_1^2}{u_1}+\frac{\xi_1}{u_1^2}\right)\frac{\partial}{\partial \xi_1}]=-\frac{a_0t_1}{u_1^2}\frac{\partial}{\partial \xi_1}\wedge \frac{\partial}{\partial u_1}-a_0t_2\left(\frac{2\xi_1}{u_1}+\frac{1}{u_1^2} \right)\frac{\partial}{\partial \xi_1}\wedge \frac{\partial}{\partial u_1}\equiv 0,\,\,\,\,\,t_1,t_2\in \mathbb{C}.
\end{align*}
Here $a\equiv b$ means that $a$ and $b$ represent the same cohomology class in $H^1(S,\wedge^2 \Theta_S)$.
Hence $\dim_\mathbb{C} \mathbb{H}^2(S,\Theta_S^\bullet)=1$, $\dim \, \mathbb{H}^1(S, \Theta_S^\bullet)=1+2=3$ by Lemma $\ref{a25}$, and
\begin{align}\label{a55}
\mathbb{H}^1(S,\Theta_S^\bullet)\cong \left\langle \frac{\partial}{\partial \xi}\wedge \frac{\partial}{\partial u}  \right\rangle \oplus \left\langle \frac{1}{u_1}\frac{\partial}{\partial u_1}, -\left(\frac{\xi_1^2}{u_1}+\frac{\xi_1}{u_1^2}\right)\frac{\partial}{\partial \xi_1} \right\rangle
\end{align}

We summarize our computation in Table $\ref{ruled}$.

\subsection{The case of $S=A_{-1}$}\

From $(\ref{a22})$, let $(S,\Lambda_0=0)$ be a Poisson ruled surface over an elliptic curve $X$. Then from $(\ref{a22})$ and Lemma $\ref{a25}$, we have $\mathbb{H}^0(S,\Theta_S^\bullet)\cong H^0(S,\Theta_S)$, $\mathbb{H}^1(S,\Theta_S^\bullet)\cong H^1(S,\Theta_S)$ and $\mathbb{H}^2(S,\Theta_S^\bullet)=0$. We summarize our computation in Table $\ref{ruled}$.

\section{Poisson deformations of ruled surfaces over an elliptic curve}\label{section4}

In this section, we determine obstructedness or unobstructedness of Poisson deformations for ruled surfaces over an elliptic curve $X$ for any Poisson structure described in section $\ref{section2}$. By extending the methods in \cite{Suw69}, in the case of  an unobstructed Poisson ruled surface $(S,\Lambda_0)$ over $X$, we show the unobstructedness in Poisson deformations by constructing a Poisson analytic family $(\mathcal{S},\Lambda,B,\pi)$ such that the associated Poisson Kodaira-Spencer map is an isomorphism at the distinguished point. Before proceeding our discussions, we recall the following lemma from \cite{Kim19}.
\begin{lemma}\label{a50}
Let $\mathcal{U}=(U \times \mathbb{P}_\mathbb{C}^1,U_1 \times \mathbb{P}_\mathbb{C}^1)$ be the Stein open covering of a ruled surface $S$ over an elliptic curve as in section $\ref{section1}$. Then $(S,\Lambda_0)$ is obstructed in Poisson deformations if for some $a , b$ where $a\in H^0(S,\wedge^2 \Theta_{S})$, and $b\in C^1(\mathcal{U},\Theta_{S})$ which defines an element in $ker(H^1(S, \Theta_{S})\xrightarrow{[\Lambda_0,-]} H^1(S,\wedge^2 \Theta_{S}))$, under the following map
\begin{align*}
[-,-]:H^0(S, \wedge^2 \Theta_{S})\times H^1(S, \Theta_{S})\to H^1(S, \wedge^2 \Theta_{S})
\end{align*}
$[a,b]\in H^1(S,\wedge^2 \Theta_{S})$ is not in the image of $ H^1(S, \Theta_{S})\xrightarrow{[\Lambda_0,-]} H^1(S,\wedge^2 \Theta_{S}) $.
\end{lemma}

We note that a group of automorphisms of $\mathbb{C}\times H$, where $H=\{\tau\in \mathbb{C}|\text{Im}\, \tau>0\}$, defined as
\begin{align*}
G=\{g_{nm}:(u,\tau )\mapsto (u+m\tau+n,\tau)|m,n\in\mathbb{Z}\}
\end{align*}
defines a complex analytic family $\sigma:\mathbb{C}\times H/G\to H$ such that $\sigma^{-1}(\omega)$ is an elliptic curve $X=\mathbb{C}/G_\omega$ as in section $\ref{section1}$, and the Kodaira-Spencer map $T_\omega H\to H^1(X,\Theta_X)\cong \mathbb{C}$ is an isomorphism at $\tau=\omega$. We note that $H^1(X,\Theta_X)\cong \langle \frac{1}{u_1}\frac{\partial}{\partial u_1}\rangle$ for $\frac{1}{u_1}\frac{\partial}{\partial u_1}\in C^1(\mathscr{U},\Theta_X)$, where $\mathscr{U}=\{U=X-p,U_1\}$ in section $\ref{section1}$. For any $(u,\tau)\in \mathbb{C}\times H$, we denote by $[(u,\tau)]$ the corresponding element of $\mathbb{C}\times H/G$. Take a point $[(p,\omega)]\in \mathbb{C}\times H/G$ and let $(u_1,\tau)$ be a local neighborhood of $[(p,\omega)]$. We set $\mathcal{V}:=\mathbb{C}\times H/G-[(p,\omega)]$ and $\mathcal{V}_1:=\{(u_1,\tau)||u_1|<\epsilon,|\tau-\omega|<\epsilon\}$. 
Then $U=X-p$ and $U_1$ in section $\ref{section1}$ satisfy $U=\mathcal{V}\cap X$ and $U_1=\mathcal{V}_1\cap X$. We will keep these notations in the following.

\subsection{The case of $S=S_0=X\times \mathbb{P}_\mathbb{C}^1$}\

We remark that in \cite{Kim19}, we already showed that $(S,\Lambda_0)$ is unobstructed in Poisson deformations for $\Lambda_0 \ne 0$, and obstructed in Poisson deformations for $\Lambda_0=0$. In \cite{Kim19}, we showed the unobstructedness by solving the integrability condition $L \alpha(t)+\frac{1}{2}[\alpha(t),\alpha(t)]=0$, where $L=\bar{\partial}+[\Lambda_0.-]$, and $\alpha(t)\in H^0(S,\wedge^2\Theta_S)\oplus H^1(S,\Theta_S)$. In this subsection, we show the unobstructedness of $(S,\Lambda_0\ne 0)$ in Poisson deformations by explicitly constructing a Poisson analytic family of deformations of $(S,\Lambda_0)$ such that the associated Poisson Kodaira-Spencer map is an isomorphism at the distinguished point. Let $(S,\Lambda_0=(A+B\xi+C\xi^2)\frac{\partial}{\partial\xi}\wedge \frac{\partial}{\partial u})$ be Poisson ruled surface over $X$ with $(A,B,C)\ne 0$ as in $(\ref{a6})$.

We set $\pi:(\mathcal{S},\Lambda)\to H\times \mathbb{C}^2$ with the coordinates $(\tau,t,t_1)\in H\times \mathbb{C}^2$ in the following way:
\begin{align*}
(\mathcal{S},\Lambda)=(\mathcal{V}\times \mathbb{P}_\mathbb{C}^1\times \mathbb{C}^2, &(A+t_1F_1+(B+t_1F_2)\xi+(C+t_1F_3)\xi^2)\frac{\partial}{\partial \xi}\wedge \frac{\partial}{\partial u})\\
&\bigcup \,\,\, (\mathcal{V}_1\times \mathbb{P}_\mathbb{C}^1\times \mathbb{C}^2,(A+t_1F_1+(B+t_1F_2)\xi_1+(C+t_1F_3)\xi_1^2)\frac{\partial}{\partial \xi_1}\wedge \frac{\partial}{\partial u_1}
)
\end{align*}
 where $([(u,\tau)],\xi,t,t_1)\in \mathcal{V}\times \mathbb{P}_\mathbb{C}^1\times \mathbb{C}^2$ and $((u_1,\tau),\xi_1,t,t_1)\in \mathcal{V}_1\times \mathbb{P}_\mathbb{C}^1\times \mathbb{C}^2$ are identified if and only if
\begin{align*}
[(u,\tau)]=(p+u_1,\tau),\,\,\,\,\,\,\xi=\frac{(1+\frac{(B+t_1F_2)t}{u_1})\xi_1+\frac{(A+t_1F_1)t}{u_1}}{-\frac{(C+t_1F_3)t}{u_1}\xi_1+1}=\frac{(1+\frac{B't}{u_1})\xi_1+\frac{A't}{u_1}}{-\frac{C't}{u_1}\xi_1+1}=\frac{(u_1+B't)\xi_1+A't}{-C't\xi_1+u_1}\\
\end{align*}
where $A':=A+t_1F_1,\,\,B':=B+t_1F_2,\,\,C':=C+t_1F_3$ and $(F_1,F_2,F_3)\in \mathbb{C}^3$ in $(\ref{a35})$. We show that $\Lambda$ is well-defined on $S$, i.e. $\Lambda\in H^0(\mathcal{S},\wedge^2 \Theta_{\mathcal{S}/H\times \mathbb{C}^2})$.  We note that
\begin{align*}
\frac{\partial}{\partial \xi}\wedge \frac{\partial}{\partial u}&=\frac{(-C't\xi_1+u_1)^2}{u_1^2 +B'tu_1+A'C't^2}\frac{\partial}{\partial \xi_1}\wedge\frac{\partial}{\partial u_1}
\end{align*}
Then we have
\begin{align*}
&(A+t_1F_1+(B+t_1F_2)\xi+(C+t_1F_3)\xi^2)\frac{\partial}{\partial \xi}\wedge \frac{\partial}{\partial u}=(A'+B'\xi+C'\xi^2)\frac{\partial}{\partial \xi}\wedge \frac{\partial}{\partial u}\\
&=\left(A'+B'\left(\frac{(u_1+B't)\xi_1+A't}{-C't\xi_1+u_1} \right)+C'\left(\frac{(u_1+B't)\xi_1+A't}{-C't\xi_1+u_1}\right)^2   \right)\frac{(-C't\xi_1+u_1)^2}{u_1^2 +B'tu_1+A'C't^2}\frac{\partial}{\partial \xi_1}\wedge \frac{\partial}{\partial u_1}\\
&=(A'+B'\xi_1+C'\xi_1^2)\frac{\partial}{\partial \xi_1}\wedge \frac{\partial}{\partial u_1}=(A+t_1F_1+(B+t_1F_2)\xi_1+(C+t_1F_3)\xi^2)\frac{\partial}{\partial \xi}\wedge \frac{\partial}{\partial u_1}
\end{align*}
Hence $\Lambda$ is well-defined. Then from $(\ref{a35})$, we see that the Poisson Kodaira-Spencer map 
\begin{align*}
T_{(\omega,0,0)}( H\times \mathbb{C}^2 )\to \mathbb{H}^1(S,\Theta_S^\bullet)
\end{align*}
is an isomorphism at $(\tau,t,t_1)=(\omega,0,0)$ . Hence $(S,\Lambda_0)$ is unobstructed in Poisson deformations. We summarize Poisson deformations of $S=X\times \mathbb{P}_\mathbb{C}^1$ in Table $\ref{ruled}$.

\subsection{The case of $S\in \mathcal{S}_0, S\ne S_0$}\

Let $(S, \Lambda_0=A\xi\frac{\partial}{\partial \xi}\wedge \frac{\partial}{\partial u})$ be a Poisson ruled surface over $X$ as in $(\ref{a36})$. We show that $(S,\Lambda_0)$ is unobstructed in Poisson deformations by constructing a Poisson analytic family of deformations of $(S,\Lambda_0)$ such that the Poisson Kodaira-Spencer map is an isomorphism at the distinguished point. We set $\pi:(\mathcal{S},\Lambda)\to H\times \mathbb{C}^2$ with the coordinates $(\tau, t,t_1)\in H\times \mathbb{C}^2$ in the following way:
\begin{align*}
(\mathcal{S},\Lambda)=(\mathcal{V}\times \mathbb{P}_\mathbb{C}^1\times \mathbb{C}^2, (A+t_1)\xi\frac{\partial}{\partial \xi}\wedge \frac{\partial}{\partial u})\,\,\,\bigcup \,\,\, (\mathcal{V}_1\times \mathbb{P}_\mathbb{C}^1\times \mathbb{C}^2, (A+t_1)\xi_1\frac{\partial}{\partial \xi_1}\wedge \frac{\partial}{\partial u_1})
\end{align*}
 where $([(u,\tau)],\xi,t,t_1)\in \mathcal{V}\times \mathbb{P}_\mathbb{C}^1\times \mathbb{C}^2$ and $((u_1,\tau),\xi_1,t,t_1)\in \mathcal{V}_1\times \mathbb{P}_\mathbb{C}^1\times \mathbb{C}^2$ are identified if and only if
\begin{align*}
[(u,\tau)]=(p+u_1,\tau),\,\,\,\,\,\,\xi=e^{\frac{t_0+t}{u_1}}\xi_1\\
\end{align*}
It is clear that $\Lambda$ is well-defined on $\mathcal{S}$. i.e. $\Lambda\in H^0(\mathcal{S},\wedge^2 \Theta_{\mathcal{S}/H\times \mathbb{C}^2})$. Then from $(\ref{a37})$, we see that the Poisson Kodaira-Spencer map  
\begin{align*}
T_{(\omega,0,0)}( H\times \mathbb{C}^2 )\to \mathbb{H}^1(S,\Theta_S^\bullet)
\end{align*}
is an isomorphism at $(\tau,t,t_1)=(\omega,0,0)$. Hence $(S,\Lambda_0)$ is unobstructed in Poisson deformations. We summarize Poisson deformations of $(S,\Lambda_0)$ in Table $\ref{ruled}$.

\subsection{The case of $S=S_n(n\geq 1)$}\

Let $(S,\Lambda_0=a_0\xi\frac{\partial}{\partial \xi}\wedge\frac{\partial}{\partial u}+A(u)\xi^2\frac{\partial}{\partial \xi}\wedge \frac{\partial}{\partial u})$ be a Poisson ruled surface over $X$ as in $(\ref{a15})$. We will show that if $a_0\ne 0$, then $(S,\Lambda_0)$ is unobstructed in Poisson deformations, and if $a_0=0$, then $(S,\Lambda_0)$ is obstructed in Poisson deformations.

Assume that $a_0\ne 0$. We set $\pi:(\mathcal{S},\Lambda)\to H\times \mathbb{C}$ with the coordinates $(\tau,t)\in H\times \mathbb{C}$ in the following way:
\begin{align*}
(\mathcal{S},\Lambda)=(\mathcal{V}\times \mathbb{P}_\mathbb{C}^1\times \mathbb{C}, \left( (a_0+t)\xi+ A(u;\tau) \xi^2\right)\frac{\partial}{\partial \xi}\wedge \frac{\partial}{\partial u})\,\,\,\bigcup \,\,\, (\mathcal{V}_1\times \mathbb{P}_\mathbb{C}^1\times \mathbb{C}, \left( (a_0+t)\xi_1+ A(u;\tau)u_1^n \xi_1^2\right)\frac{\partial}{\partial \xi_1}\wedge \frac{\partial}{\partial u_1})
\end{align*}
 where $([(u,\tau)],\xi,t)\in \mathcal{V}\times \mathbb{P}_\mathbb{C}^1\times \mathbb{C}$ and $((u_1,\tau),\xi_1,t)\in \mathcal{V}_1\times \mathbb{P}_\mathbb{C}^1\times \mathbb{C}$ are identified if and only if
\begin{align*}
[(u,\tau)]=(p+u_1,\tau),\,\,\,\,\,\,\xi=u_1^n \xi_1,
\end{align*}
and $A(u:\tau):=c_0+c_1\wp(u-p;\tau)+c_2\wp'(u-p;\tau)+\cdots+c_{n-1}\wp^{(n-2)}(u-p;\tau)$. Here $\wp(u;\tau)$ is the Weierstrass $\wp$-function with periods $(1,\tau)$. We note that $A(u)=A(u;\omega)$.
It is clear that $\Lambda$ is well-defined on $\mathcal{S}$. i.e. $\Lambda\in H^0(\mathcal{S},\wedge^2 \Theta_{\mathcal{S}/H\times \mathbb{C}})$. Then from $(\ref{a39})$, we see that the Poisson Kodaira-Spencer map  
\begin{align*}
T_{(\omega,0)}( H\times \mathbb{C} )\to \mathbb{H}^1(S,\Theta_S^\bullet)
\end{align*}
is an isomorphism at $(\tau,t)=(\omega,0)$. Hence $(S,\Lambda_0)$ is unobstructed in Poisson deformations.

Next assume that $a_0=0, n=1$ (i.e. $S=S_1$) and $A(u)=c_0\ne 0$. We set $\pi:(\mathcal{S},\Lambda)\to H\times \mathbb{C}$ with the coordinates $(\tau,t)\in H\times \mathbb{C}$ in the following way:
\begin{align*}
(\mathcal{S},\Lambda)=(\mathcal{V}\times \mathbb{P}_\mathbb{C}^1\times \mathbb{C}, \left( t\xi+ c_0\xi^2\right)\frac{\partial}{\partial \xi}\wedge \frac{\partial}{\partial u})\,\,\,\bigcup \,\,\, (\mathcal{V}_1\times \mathbb{P}_\mathbb{C}^1\times \mathbb{C}, \left( t\xi_1+ c_0u_1 \xi_1^2\right)\frac{\partial}{\partial \xi_1}\wedge \frac{\partial}{\partial u_1})
\end{align*}
 where $([(u,\tau)],\xi,t)\in \mathcal{V}\times \mathbb{P}_\mathbb{C}^1\times \mathbb{C}$ and $((u_1,\tau),\xi_1,t)\in \mathcal{V}_1\times \mathbb{P}_\mathbb{C}^1\times \mathbb{C}$ are identified if and only if
\begin{align*}
[(u,\tau)]=(p+u_1,\tau),\,\,\,\,\,\,\xi=u_1 \xi_1,
\end{align*}
Then from $(\ref{a90})$, we see that the Poisson Kodaira-Spencer map  
\begin{align*}
T_{(\omega,0)}( H\times \mathbb{C} )\to \mathbb{H}^1(S_1,\Theta_{S_1}^\bullet)
\end{align*}
is an isomorphism at $(\tau,t)=(\omega,0)$. Hence $(S=S_1,\Lambda_0)$ is unobstructed in Poisson deformations.

On the other hand, assume that $\Lambda_0=0$ (i.e $a_0=0,A(u)=0$) for $n \geq 1$ or $a_0=0,A(u)\ne 0$ for $n\geq 2$. Then from $(\ref{a51})$, consider
\begin{align*}
[\xi_1\frac{\partial}{\partial \xi_1}\wedge \frac{\partial}{\partial u_1},\left( \frac{t_1}{u_1^{n+1}}+\frac{t_2}{u_1}+\frac{t_3}{u_1^2}+\cdots +\frac{t_n}{u_1^{n-1}}\right)\frac{\partial}{\partial \xi_1}]\equiv  -\left( \frac{t_1}{u_1^{n+1}}+\frac{t_2}{u_1}+\frac{t_3}{u_1^2}+\cdots +\frac{t_n}{u_1^{n-1}}\right)\frac{\partial}{\partial \xi_1}\wedge \frac{\partial}{\partial u_1}
\end{align*}
where $(t_1,...,t_n)\ne 0\in \mathbb{C}^n$, which is not in the image $H^1(S,\Theta_S)\xrightarrow{[\Lambda_0,-]} H^1(S,\wedge^2 \Theta_S)$ since $a_0=0$. Then by Lemma $\ref{a50}$, $(S,\Lambda_0)$ is obstructed in Poisson deformations. We summarize Poisson deformations of $(S,\Lambda_0)$ in Table $\ref{ruled}$.

\subsection{The case of $S=A_0$}\

Let $(S,\Lambda_0=a_0\frac{\partial}{\partial \xi}\wedge \frac{\partial}{\partial u})$ be a Poisson ruled surface over $X$ as in $(\ref{a53})$. We show that $(S,\Lambda_0)$ is unobstructed in Poisson deformations by constructing a Poisson analytic family of deformations of $(S,\Lambda_0)$ such that the Poisson Kodaira-Spencer map is an isomorphism at the distinguished point. We set $\pi:(\mathcal{S},\Lambda)\to H\times \mathbb{C}^2$ with the coordinates $(\tau,t,t_1)\in H\times \mathbb{C}^2$ in the following way:
\begin{align*}
(\mathcal{S},\Lambda)=(\mathcal{V}\times \mathbb{P}_\mathbb{C}^1\times \mathbb{C}^2, (a_0+t_1)(1-t\xi^2)\frac{\partial}{\partial \xi}\wedge \frac{\partial}{\partial u})\,\,\,\bigcup \,\,\, (\mathcal{V}_1\times \mathbb{P}_\mathbb{C}^1\times \mathbb{C}^2, (a_0+t_1)(1-t\xi_1^2)\frac{\partial}{\partial \xi_1}\wedge \frac{\partial}{\partial u_1})
\end{align*}
 where $([(u,\tau)],\xi,t,t_1)\in \mathcal{V}\times \mathbb{P}_\mathbb{C}^1\times \mathbb{C}^2$ and $((u_1,\tau),\xi_1,t,t_1)\in \mathcal{V}_1\times \mathbb{P}_\mathbb{C}^1\times \mathbb{C}^2$ are identified if and only if
\begin{align*}
[(u,\tau)]=(p+u_1,\tau),\,\,\,\,\,\,\xi=\frac{\xi_1+\frac{1}{u_1}}{\frac{t}{u_1}\xi_1+1}=\frac{u_1\xi_1+1}{t\xi_1+u_1}\\
\end{align*}
We show that $\Lambda$ is well-defined on $\mathcal{S}$. i.e. $\Lambda\in H^0(\mathcal{S},\wedge^2 \Theta_{\mathcal{S}/H\times \mathbb{C}^2})$. We note that
\begin{align*}
\frac{\partial}{\partial \xi}\wedge \frac{\partial}{\partial u}&=\frac{(t\xi_1+u_1)^2}{u_1^2-t}\frac{\partial}{\partial \xi_1}\wedge \frac{\partial}{\partial u_1},\\
1-t\xi^2&=1-t\left(\frac{u_1\xi_1+1}{t\xi_1+u_1} \right)^2=\frac{(1-t\xi_1^2)(u_1^2-t)}{(t\xi_1+u_1)^2}
\end{align*}
Then we have
\begin{align*}
(a_0+t_1)(1-t\xi^2)\frac{\partial}{\partial \xi}\wedge \frac{\partial}{\partial u}=(a_0+t_1)\left(\frac{(1-t\xi_1^2)(u_1^2-t)}{(t\xi_1+u_1)^2}\right)\frac{(t\xi_1+u_1)^2}{u_1^2-t}\frac{\partial}{\partial \xi_1}\wedge \frac{\partial}{\partial u_1}
=(a_0+t_1)(1-t\xi_1^2)\frac{\partial}{\partial \xi_1}\wedge \frac{\partial}{\partial u_1}
\end{align*}
Hence $\Lambda$ is well-defined.
Then from $(\ref{a55})$, we see that the Poisson Kodaira-Spencer map  
\begin{align*}
T_{(\omega,0,0)}( H\times \mathbb{C}^2 )\to \mathbb{H}^1(S,\Theta_S^\bullet)
\end{align*}
is an isomorphism at $(\tau,t,t_1)=(\omega,0,0)$. Hence $(S,\Lambda_0)$ is unobstructed in Poisson deformations. We summarize Poisson deformations of $(S=A_0,\Lambda_0)$ in Table $\ref{ruled}$.

\subsection{The case of $S=A_{-1}$}\

Since there is no nontrivial Poisson structure on $S$ and $S$ is unobstructed in complex deformations, $S$ is unobstructed in Poisson deformations. We summarize Poisson deformations of $(S,\Lambda_0=0)$ in Table $\ref{ruled}$.

\begin{table}
\begin{center}
\begin{tabular}{| c | c | c | c | c | c | c | } \hline
Type of $S$ &  Poisson structure $\Lambda_0$ & $\dim_\mathbb{C}  \mathbb{H}^0(S, \Theta_S^\bullet)$ & $\dim_\mathbb{C} \mathbb{H}^1(S,\Theta_S^\bullet)$ &$\dim_\mathbb{C} \mathbb{H}^2(S,\Theta_S^\bullet)$ & Poisson deformations  \\ \hline
$S_0=X\times \mathbb{P}_\mathbb{C}^1$ & $0$ & $4$ & $7$  & $3$  &  obstructed \\ \hline
$S_0=X\times \mathbb{P}_\mathbb{C}^1$ & $(A,B,C)\ne 0$ in $(\ref{a6})$ & $2$ & $3$  & $1$   & unobstructed \\ \hline
$S\in \mathcal{S}_0, S\ne S_0$ & any Poisson structure  & $2$ & $3$  & $1$   & unobstructed \\ \hline
$S_n (n\geq 1)$ & $ 0 $ & $n+1$ & $2n+2$  & $n+1$  & obstructed \\ \hline
$S_1$ & $a_0=0, c_0 \ne 0$ in $(\ref{a15})$ & $1$ & $2$  & $1$  & unobstructed \\ \hline
$S_n (n\geq 2)$ & $a_0=0, A(u)\ne 0$ in $(\ref{a15})$  & $n$ &  $2n$ & $n$  & obstructed \\ \hline
$S_n (n\geq 1)$ & $a_0\ne 0$ in $(\ref{a15})$  & $1$ & $2$ & $1$ & unobstructed \\ \hline
$A_0$ & any Poisson structure  & $2$ & $3$  & $1$ & unobstructed  \\ \hline
$A_{-1}$ & $0$  & $1$ & $1$  & $0$ & unobstructed \\ \hline
\end{tabular}
\end{center}
\caption{obstructed and unobstructedness of Poisson ruled surfaces $(S,\Lambda_0)$ over an elliptic curve $X$} \label{ruled}
\end{table}

\bibliographystyle{amsalpha}
\bibliography{References-Rev9} 

\end{document}